\documentclass[conference]{IEEEtran}

\usepackage{amsmath,amssymb,amsthm}
\usepackage{graphicx}
\usepackage{bm,color}
\usepackage{url}
\usepackage{cite}
\usepackage{textcomp}
\usepackage{slashbox}
\usepackage{bigints}

\newtheorem{theorem}{Theorem}

\newtheorem{proposition}{Proposition}
\theoremstyle{definition}
\newtheorem{remark}{Remark}

\newcommand{\st}{\textrm{s.t.}}

\newcommand{\diff}{\mathop{}\mathrm{d}}

\newcommand{\cT}{{\cal T}}

\makeatletter
\def\hlinew#1{%
  \noalign{\ifnum0=`}\fi\hrule \@height #1 \futurelet
   \reserved@a\@xhline}

\DeclareMathOperator*{\argmin}{arg\,min}

\IEEEoverridecommandlockouts

\begin{document}
\title{Robust Optimal Power Flow with Wind Integration Using Conditional Value-at-Risk}

\author{\IEEEauthorblockN{Yu Zhang and Georgios B. Giannakis}
\authorblockA{Dept. of ECE and DTC, University of Minnesota, Minneapolis, USA \\
Emails: \{zhan1220, georgios\}@umn.edu}
\thanks{This work was supported by the
NSF ECCS grant 1202135, and University of Minnesota
Institute of Renewable Energy and the Environment (IREE) grant RL-0010-13.}}

\maketitle

\begin{abstract}
Integrating renewable energy into the power grid requires intelligent risk-aware dispatch accounting for the stochastic availability of renewables.
Toward achieving this goal, a robust DC optimal flow problem is developed in the present paper for power systems with a high penetration of wind energy.
The optimal dispatch is obtained as the solution to a convex program with a suitable regularizer, which is able to mitigate the potentially high
risk of inadequate wind power. The regularizer is constructed based on the energy transaction cost using conditional value-at-risk (CVaR).
Bypassing the prohibitive high-dimensional integral, the distribution-free sample average approximation method is efficiently utilized
for solving the resulting optimization problem. Case studies are reported to corroborate the efficacy of the novel model and approach tested on the
IEEE $30$-bus benchmark system with real operation data from seven wind farms.
\end{abstract}

\section{Introduction}
With the appealing attributes of being environment-friendly and price-competitive over conventional power generation,
clean renewable sources of energy, such as wind, solar, hydro, and geothermal energy, have been developing rapidly over the last few decades.
Growing at an annual rate of $20$\%, wind power generation had $282.5$ GW worldwide installed capacity by the end of $2012$~\cite{GWEC}.
The U.S. Department of Energy set a goal of using wind energy to generate $20$\% of the nation's electricity demand by $2030$~\cite{DOE08}.

Aligned to the goal of boosting the penetration of renewable energy sources in future smart grids,
energy management with renewables, including economic dispatch (ED), unit commitment (UC),
and optimal power flow (OPF), have been extensively investigated recently.
Generally, two types of strategies have been developed to address the key challenge of
dealing with the supply-demand balance, which is induced by the stochastic availability
and intermittency of renewables. Early works aim at maintaining balance by introducing
committed renewable energy. ED penalizing overestimation and underestimation
of wind power is investigated in~\cite{HetzerYB08}. Worst-case robust distributed dispatch
with demand side management is proposed for grid-connected microgrids
with distributed energy resources in~\cite{YuNGGG-TSE13}.
Its solution though can be very sensitive to the accuracy of the wind power forecast.
For the second type, supply-demand imbalance is allowed up to a certain extent by limiting the loss-of-load probability (LOLP).
Leveraging the scenario approximation technique, risk-limiting ED and DC-OPF tasks with correlated wind farms
have been developed recently in~\cite{YuNGGG-ISGT13} and~\cite{YuNKGG-DSP13}, respectively.
A multi-stage stochastic control problem is pursued for risk-limiting dispatch of wind power in~\cite{RajagopalACC12}.
Chance-constrained two-stage stochastic program is formulated in~\cite{FangGW12} for UC with uncertain wind power output;
see also~\cite{BiChHa13} and~\cite{SjGaTo12} for advances in chance-constrained OPF.
However, the applied conic optimization technique therein relies on Gaussianity assumptions for the wind power generation.

Additional limitations are present in the aforementioned works. For example, worst-case renewable
energy generation is unlikely when it comes to real time operation~\cite{YuNGGG-TSE13}.
The chance-constrained problems are typically non-convex for general probability distributions.
Leveraging the scenario sampling, the relaxed convex problems can be solved efficiently.
However, in certain scenarios, this technique turns out to be too conservative for scheduling
the delivered renewables~\cite{YuNGGG-ISGT13}.

This paper deals with robust DC OPF for a smart grid with high penetration of wind power.
Instead of a chance-constrained formulation, an optimization problem is introduced with an
appropriate regularizer that plays an instrumental role for mitigating the high risk of inadequate wind power.
The regularizer is constructed based on the energy transaction cost using the \emph{conditional value-at-risk} (CVaR).
This ``smart'' CVaR-based regularizer turns out to be capable of utilizing renewables intelligently with limited risk.
The resulting optimization problem, which aims at minimizing jointly the generation as well as the energy transaction costs,
is provably convex thanks to the CVaR.
To bypass the prohibitive high-dimensional integral present in the regularizer,
the sample average approximation (SAA) is utilized to obtain an efficient distribution-free approach.
Numerical tests are performed on the IEEE $30$-bus benchmark system to corroborate the effectiveness of the novel
model and approach using real wind farm data~\cite{PSTCA}, \cite{MATPOWER}, \cite{kaggle}.

The remainder of the paper is organized as follows. Section~\ref{sec:CVAR} introduces VaR and CVaR.
Section~\ref{sec:Probformulation} formulates the CVaR-based DC-OPF problem along with the SAA solver.
Numerical results are reported in Section~\ref{sec:Numericalresults}, while conclusions and research directions
can be found in Section~\ref{sec:Conclusion}.

\section{Risk Measure: VaR and CVaR}~\label{sec:CVAR}
Being the most popular measures of risk, the value-at-risk (VaR) and the conditional value-at-risk (CVaR)
play a central role in risk-aware portfolio optimization~\cite{Fabozzi07}.
These risk metrics were introduced in the ground-breaking works of~\cite{RTR00} and~\cite{RTR02}. 
The redux here is useful to grasp the role of these metrics in the present context.

Let the real-valued loss function $L(\mathbf{x},\bm{\xi}) : X \times \Xi \mapsto \mathbb{R}$ denote the cost associated with the
decision variable $\mathbf{x} \in X \subset \mathbb{R}^n$; and the random vector $\bm{\xi}$ with probability density function $p(\bm{\xi})$ supported on a set $\Xi \subset \mathbb{R}^d$.
In the context of power systems, $\mathbf{x}$ can represent for instance the power schedules of conventional generators while $\bm{\xi}$ captures the sources of uncertainty due to e.g., renewables, forecasted load demand, and locational marginal prices (LMPs). The operator-concerned loss $L(\mathbf{x},\bm{\xi})$ represents the cost, which depends on both $\mathbf{x}$ and $\bm{\xi}$.
Clearly, the probability of $L(\mathbf{x},\bm{\xi})$ not exceeding a threshold $\eta$ is given by
\begin{align}
\Psi(\mathbf{x},\eta) = \int\limits_{L(\mathbf{x},\bm{\xi})\le \eta}\!\! p(\bm{\xi})\diff \bm{\xi}.
\end{align}
It can be seen that $\Psi$ is the cumulative distribution function (CDF) for the loss parameterized by $\mathbf{x}$,
which is right continuous and nondecreasing in $\eta$.
Let $\eta_\beta(\mathbf{x})$ and $\phi_\beta(\mathbf{x})$ denote respectively the $\beta$-VaR and $\beta$-CVaR values of the random loss with
a prescribed probability level $\beta \in (0,1)$. Commonly chosen values of $\beta$ are, e.g., $0.99$, $0.95$, and $0.9$.
Dependent on $\Psi$, the $\beta$-VaR and $\beta$-CVaR values are defined as
\begin{align}
\eta_\beta(\mathbf{x}) &:= \min\{\eta \in\mathbb{R}~|~\Psi(\mathbf{x},\eta) \ge \beta\}\label{eq:VaR} \\
\phi_\beta(\mathbf{x}) &:= \frac{1}{(1-\beta)}\!\int\limits_{L(\mathbf{x},\bm{\xi})\ge \eta_\beta(\mathbf{x})}\!\! L(\mathbf{x},\bm{\xi})p(\bm{\xi})\diff \bm{\xi}.\label{eq:CVaR}
\end{align}
Since $\Psi$ is non-decreasing in $\eta$, $\eta_\beta(\mathbf{x})$
comes out as the left endpoint of the nonempty interval consisting of the solution $\eta$ satisfying $\Psi(\mathbf{x},\eta) = \beta$.
Hence, $\phi_\beta(\mathbf{x})$ is the conditional expectation of $L(\mathbf{x},\bm{\xi})$ to be greater than or equal to $\eta_\beta(\mathbf{x})$.

The characterization of $\eta_\beta(\mathbf{x})$ and $\phi_\beta(\mathbf{x})$ lies in the optimization of a key constructed function
\begin{align}
F_{\beta}(\mathbf{x},\eta) = \eta +\frac{1}{1-\beta}\int_\mathbf{\bm{\xi}\in \Xi}\left[L(\mathbf{x},\bm{\xi})-\eta\right]^{+}p(\bm{\xi})\diff \bm{\xi} \label{eq:F}
\end{align}
where $[a]^{+}:=\max\{a,0\}$ is the projection operator.
The crucial features of $F_{\beta}$ relating $\eta_\beta(\mathbf{x})$ with $\phi_\beta(\mathbf{x})$ are summarized in the following theorem.
\begin{theorem}[\cite{RTR00}, pp.~24--26]\label{Them:RTR}
Function $F_{\beta}(\mathbf{x},\eta)$ is convex and continuously differentiable in $\eta$. Furthermore,
$F_{\beta}(\mathbf{x},\eta)$ is convex with respect to $(\mathbf{x},\eta)$ while $\phi_\beta(\mathbf{x})$ is convex in $\mathbf{x}$, provided that
$L(\mathbf{x},\bm{\xi})$ is convex in $\mathbf{x}$. The relationships among $F_{\beta}(\mathbf{x},\eta)$, $\eta_\beta(\mathbf{x})$, and $\phi_\beta(\mathbf{x})$ are given as follows
\begin{align}
\phi_\beta(\mathbf{x}) &= \min_{\eta \in \mathbb{R}}\,F_{\beta}(\mathbf{x},\eta) \\
\eta_\beta(\mathbf{x}) &= \lfloor\argmin_{\eta \in \mathbb{R}}\,F_{\beta}(\mathbf{x},\eta)\rfloor \\
\min_{\mathbf{x}\in X}\phi_\beta(\mathbf{x}) &= \min_{(\mathbf{x},\eta)\in X\times \mathbb{R}}F_{\beta}(\mathbf{x},\eta)
\end{align}
where $\lfloor\Pi\rfloor$ denotes the left endpoint of the interval set $\Pi$.
\end{theorem}

It is important to appreciate the claim in Theorem~\ref{Them:RTR} regarding the undesirable characteristics of $\beta$-VaR, 
namely non-subadditivity and non-convexity.
Theorem~\ref{Them:RTR} asserts that minimizing the convex $\beta$-CVaR $\phi_\beta(\mathbf{x})$ is equivalent to minimizing $F_{\beta}(\mathbf{x},\eta)$, which is not only convex,
but also easier to approximate.
A straightforward and easily implementable approximation of the expectation function $F_{\beta}$ is its empirical estimate using $N_s$ Monte Carlo samples $\{\bm{\xi}_s\}_{s=1}^{N_s}$,
namely
\begin{align}
\hat{F}_{\beta}(\mathbf{x},\eta) = \eta +\frac{1}{N_s(1-\beta)}\sum_{s=1}^{N_s}\left[L(\mathbf{x},\bm{\xi}_s)-\eta\right]^{+}.\label{eq:approxF}
\end{align}
Clearly, the sample average approximation method is distribution free, and
the law of large numbers asserts $\hat{F}_{\beta}$ as a good approximation of ${F}_{\beta}$ for $N_s$ large enough.
Furthermore, $\hat{F}_{\beta}(\mathbf{x},\eta)$ is convex with respect to $(\mathbf{x},\eta)$ when $L(\mathbf{x},\bm{\xi}_s)$ is convex in $\mathbf{x}$.
The non-differentiability  due to the projection operator can be readily overcome by leveraging the epigraph form of $\hat{F}$, which will be shown explicitly in Section~\ref{sec:probStat}.

Leveraging CVaR, a robust OPF problem will be formulated next by considering the transaction cost induced by wind power shortage.


\section{Robust Optimal Power Flow Formulation}\label{sec:Probformulation}

Consider a power grid with $M$ buses. Let $\mathbf{p}_G :=[p_{G_1},\ldots,p_{G_M}]^\cT$ denote a vector collecting 
the conventional power outputs of the thermal generators,
and $\mathbf{p}_D :=[p_{D_1},\ldots,p_{D_M}]^\cT$ the load demand, where $(\cdot)^\cT$ denotes transposition.
Furthermore, if a renewable energy facility (e.g., a wind farm) is located at bus $m$ as well, two quantities will be associated with it: the actual wind power generation $w_m$,
and the power $p_{W_m}$ scheduled to be injected to bus $m$. Note that the former is a random variable, whereas the latter is a decision variable.
For notational simplicity, define further two $M$-dimensional vectors
$\mathbf{w} :=[w_1,\ldots,w_M]^\cT$, and $\mathbf{p}_W := [p_{W_1},\ldots,p_{W_M}]^\cT$.
Clearly, if no generator, load, or wind farm is attached to bus $m$, the $m$th entry of $\mathbf{p}_G$, $\mathbf{p}_D$, 
or, $\mathbf{w}$ and $\mathbf{p}_W$ is set to zero.

\subsection{DC Flow Model}

For a power transmission network with a total of $N$ lines, let $x_n$ denote the reactance associated with the $n$th line.
Define further a diagonal matrix $\mathbf{D} := \mathrm{diag}\left(x_1^{-1},\ldots,x_N^{-1}\right) \in \mathbb{R}^{N \times N}$, and the branch-bus incidence
matrix $\mathbf{A} \in \mathbb{R}^{N \times M}$, such that if its $n$th row $\mathbf{a}_n^\cT$
corresponds to the branch $(i,j)$, then $[\mathbf{a}_n]_i:=+1$, $[\mathbf{a}_n]_j:=-1$, and zero elsewhere.

Consider now the DC flow model~\cite{ExpConCanBook}, and let vector $\bm\theta:=[\theta_1,\ldots,\theta_M]^\cT$ collect the nodal voltage phases $\{\theta_m\}_{m=1}^M$.
Then, the power flows on all transmission lines can be expressed as $\mathbf{H}\bm\theta$ with $\mathbf{H} := \mathbf{D} \mathbf{A}$.
Physical considerations enforce a power flow limit $\mathbf{f}^{\max}$ on each transmission line, leading to the \emph{line flow constraint}
\begin{align*}
-\mathbf{f}^{\max} \preceq \mathbf{H}\bm\theta \preceq \mathbf{f}^{\max}
\end{align*}
where $\preceq$ denotes entry-wise inequality.

Furthermore, flow conservation compels zero net flow at each bus; i.e., the outgoing power flow must equal the aggregate incoming power flows.
This gives rise to the \emph{nodal balance equation} for the DC power flow model:
\begin{align}
\label{eq:node-bal}
\mathbf{p}_G + \mathbf{p}_W - \mathbf{p}_D = \mathbf{B}\bm\theta
\end{align}
where $\mathbf{B} := \mathbf{A}^\cT \mathbf{D} \mathbf{A}$ is the bus admittance matrix.
With $\mathbf{1}$ denoting the all-ones vector, it holds that $\mathbf{B} \cdot \mathbf{1} = \mathbf{0}$,
which implies that~\eqref{eq:node-bal} is invariant to nodal phase shifts.
Hence, without loss of generality, the first bus can be set to be the zero phase reference bus, i.e., $\theta_1=0$.
For simplicity, only non-dispatchable base loads will be included in $\mathbf{p}_D$.
These are fixed constants for the optimization problem that will be formulated later.
\begin{remark}[Actual versus committed wind power]
Since day-ahead power dispatch is considered in this paper, power generation schedules
should be decided prior to real time operation. Hence, the actual wind power output $\mathbf{w}$, 
which is random due to the wind speed variability, is not available at the decision making time. 
However, to maintain node balance, slack variables $\mathbf{p}_W$ will be introduced to capture 
the committed wind power injected at the corresponding buses.
These are possible to determine before the real time operation, 
together with other decision variables, namely $\mathbf{p}_G$ and $\bm{\theta}$.
\end{remark}

\subsection{CVaR-based Energy Transaction Cost}

Since the wind generation output is stochastic, it is unlikely that the scheduled power $\mathbf{p}_W$ will be equal to the actual one $\mathbf{w}$.
Thus, in order to satisfy the nodal balance~\eqref{eq:node-bal} in real time operation, either energy surplus or shortage should be included.
In the former case, the wind generation company (W-GENCO) may simply choose to curtail the excess wind power at almost no cost.
For the case of shortage, in order to accomplish the bid as promised in its signed day-ahead contract, W-GENCO has then to buy the energy shortfall
from the real time market in the form of ancillary services. Generally, wind farms attached to different buses may have different purchase prices.
This is simply because they may resort to different energy sellers, or because of the varying real-time LMPs across the grid.

Let $T_m$ denote the purchase transaction cost for the renewable energy facility associated with the $m$th bus.
Clearly, with the power shortfall being $[p_{W_m}-w_m]^{+}$ at bus $m$, the grid-wide total transaction cost is given by
$T(\mathbf{p}_W,\mathbf{w}) =\sum_{m=1}^{M} T_m\left([p_{W_m}-w_m]^{+}\right)$.
If the general loss function $L(\cdot,\cdot)$ in~\eqref{eq:F} is replaced by the transaction cost $T(\cdot,\cdot)$,
function $F_{\beta}$ related to the conditional expected transaction cost turns out to be
\begin{align}
F_{\beta}(\mathbf{p}_W,\eta) = \eta +\frac{1}{1-\beta}\mathbb{E}_{\mathbf{w}}\left[\sum_{m=1}^M T_m\left([p_{W_m}-w_m]^{+}\right)-\eta\right]^{+} \label{eq:F-pW}
\end{align}
where $\mathbb{E}[\cdot]$ denotes expectation. The following proposition sheds light on the convexity of $F_{\beta}(\mathbf{p}_W,\eta)$.
\begin{proposition}
\label{prop:convex}
If all costs $\{T_m(\cdot)\}_{m=1}^M$ are convex and non-decreasing, then $F_{\beta}(\mathbf{p}_W,\eta)$ is convex with respect to $(\mathbf{p}_W,\eta)$.
\end{proposition}
\begin{IEEEproof}
Thanks to Theorem~\ref{Them:RTR}, it suffices to show that $T(\mathbf{p}_W,\mathbf{w}) =\sum_{m=1}^{M} T_m\left([p_{W_m}-w_m]^{+}\right)$ is convex in $\mathbf{p}_W$.
Clearly, as a pointwise maximum operation, $[p_{W_m}-w_m]^{+} = \max\{p_{W_m}-w_m, 0\}$ is convex in $p_{W_m}$.
Thus, by the convexity composition rule~\cite[Sec.~3.2.4]{Boyd}, $T_m\left([p_{W_m}-w_m]^{+}\right)$ is convex
in $p_{W_m}$ whenever $T_m(\cdot)$ is convex and non-decreasing.
The claim follows immediately with the final summation operation.
\end{IEEEproof}
It is worth pointing out that cost functions $\{T_m(\cdot)\}_{m=1}^M$ typically satisfy the condition of Proposition~\ref{prop:convex}.
In the simple linear case for which $T_m([p_{W_m}-w_m]^{+}) = c_{W_m}[p_{W_m}-w_m]^{+}$, each constant $c_{W_m}\ge 0$
actually denotes the purchase price at bus $m$.
If $\mathbf{c}_W := [c_{W_1},\ldots,c_{W_M}]^{\cT}$, then $F_{\beta}$ in~\eqref{eq:F-pW} can be written as
\begin{align}
F_{\beta}(\mathbf{p}_W,\eta) = \eta +\frac{1}{1-\beta}\mathbb{E}_{\mathbf{w}}\left[\mathbf{c}_W^\cT[\mathbf{p}_{W}-\mathbf{w}]^{+}-\eta\right]^{+}. \label{eq:F-pW-linear}
\end{align}

It is now possible to formulate the robust DC-OPF task with the CVaR-based transaction cost, as in the ensuing section.

\subsection{Problem Statement}\label{sec:probStat}

Let $C_m(p_{G_m})$ be the generation cost associated with the $m$th thermal generator. Function $C_m(p_{G_m})$ is chosen convex, typically quadratic or piecewise linear.
The robust DC-OPF problem amounts to minimizing the conventional generation cost, as well as the CVaR-based transaction cost
under certain physical grid operation constraints; that is,
\begin{subequations}
\label{eq:DCOPF-all}
\begin{align}
\text{(P1)}\quad \min_{\mathbf{p}_G, \mathbf{p}_W,\bm{\theta},\eta}~& \sum_{m=1}^{M} C_m(p_{G_m}) + \mu F_{\beta}(\mathbf{p}_W,\eta)
\label{eq:DCOPF-obj}\\
\st\quad\quad & -\mathbf{f}^{\max} \preceq \mathbf{H}\bm\theta \preceq \mathbf{f}^{\max} \label{eq:DCOPF-line}\\
& \theta_1=0 \label{eq:DCOPF-ref}\\
& \mathbf{p}_G + \mathbf{p}_W - \mathbf{p}_D = \mathbf{B}\bm\theta \label{eq:DCOPF-node} \\
& \mathbf{p}_{G}^{\min} \preceq \mathbf{p}_{G} \preceq \mathbf{p}_{G}^{\max}\label{eq:DCOPF-gen}\\
& \mathbf{p}_W \succeq \mathbf{0} \label{eq:DCOPF-pW}
\end{align}
\end{subequations}
where the risk-aversion parameter $\mu>0$ controls the trade off between
the generation cost and the transaction cost,
which should be pre-determined based on the operator's concern.
Besides constraints~\eqref{eq:DCOPF-line} --~\eqref{eq:DCOPF-node},
constraints~\eqref{eq:DCOPF-gen} and~\eqref{eq:DCOPF-pW} entail the physical limits of $\mathbf{p}_G$ and $\mathbf{p}_W$, respectively,
namely $\mathbf{p}_{G}^{\min}  := [p_{G_1}^{\min},\ldots,p_{G_M}^{\min}]^\cT$ and $\mathbf{p}_{G}^{\max}  := [p_{G_1}^{\max},\ldots,p_{G_M}^{\max}]^\cT$.
Only a single scheduling period is considered here. However, \eqref{eq:DCOPF-all} can be readily extended to formulate multi-period dispatch
with time coupling constraints, e.g., ramping up/down rates and unit minimum-up/down constraints (see e.g.,~\cite{FangGW12}).

Alternatively, it is reasonable to consider a CVaR-constrained problem,
which minimizes the generation cost with a constraint to ensure that the conditional expected
transaction cost is no more than a given budget $b$. The corresponding problem formulation can be written as
\begin{subequations}
\label{eq:DCOPF2}
\begin{align}
\text{(P2)}\quad \min_{\mathbf{p}_G, \mathbf{p}_W,\bm{\theta},\eta}~& \sum_{m=1}^{M} C_m(p_{G_m}) \label{eq:DCOPF2-obj}\\
\st\quad\quad & \eqref{eq:DCOPF-line}-\eqref{eq:DCOPF-pW}\\
& F_{\beta}(\mathbf{p}_W,\eta) \le b. \label{eq:DCOPF2-CVaRconstr}
\end{align}
\end{subequations}
\begin{remark}[Interpretation as risk-limiting dispatch]
(P1) extends the standard DC OPF problem (see e.g., \cite{ChWoWa00}) to account for uncertain wind integration.
From the perspective of  nodal balance, it is desirable to inject $\{p_{W_m}\}_{m=1}^M$ as much as possible, so that
the generation cost $\sum_{m=1}^{M} C_m(p_{G_m})$ can be potentially reduced with the decreased $\{p_{G_m}\}_{m=1}^M$.
However, increasing $\{p_{W_m}\}_{m=1}^M$ will increase the CVaR-based transaction cost $F_{\beta}(\mathbf{p}_W,\eta)$
since it is non-decreasing in $\{p_{W_m}\}_{m=1}^M$ [cf.~\eqref{eq:F-pW}].
Hence, in this sense, the regularizer $F_{\beta}(\mathbf{p}_W,\eta)$ can be interpreted as a penalty to reduce
the high transaction cost due to wind power shortage.
Finally, (P1) can be also regarded as the equivalent Lagrangian form of (P2), provided that $\mu$ is the
Lagrange multiplier corresponding to the CVaR constraint~\eqref{eq:DCOPF2-CVaRconstr}.
\end{remark}

It is clear that under the condition of Proposition~\ref{prop:convex}, the objective as well as the constraints of (P1) (and also (P2)) are all convex,
which makes (P1) and (P2) also easy to solve in principle. Nevertheless, due to the high-dimensional integral present
in $F_{\beta}(\mathbf{p}_W,\eta)$ [cf.~\eqref{eq:F-pW} and~\eqref{eq:F}], an analytical solution is tough.
To this end, it is necessary to re-write the resulting problem in a form suitable for off-the-shelf solvers.

Without loss of generality, consider (P1) with the CVaR-based regularizer $F_{\beta}$ given by~\eqref{eq:F-pW-linear}.
First, as shown in~\eqref{eq:approxF}, an efficient approximation of~\eqref{eq:F-pW}
is the empirical expectation via samples $\{\mathbf{w}_s\}_{s=1}^{N_s}$, which is given by
\begin{align}
\hat{F}_{\beta}(\mathbf{p}_W,\eta) = \eta +\frac{1}{N_s(1-\beta)}\sum_{s=1}^{N_s}\left[\mathbf{c}_W^\cT[\mathbf{p}_{W}-\mathbf{w}_s]^{+}-\eta\right]^{+}.
\label{eq:F-linear-appro}
\end{align}

Next, introduce auxiliary variables $\{\mathbf{v}_s\}_{s=1}^{Ns}$ to first upper bound the inner projection
terms $\left\{[\mathbf{p}_{W}-\mathbf{w}_s]^{+}\right\}_{s=1}^{Ns}$. Then, further upper bound the resulting
terms $\left\{\left[\mathbf{c}_W^\cT\mathbf{v}_s-\eta\right]^{+}\right\}_{s=1}^{Ns}$
using another group of auxiliary variables $\{u_s\}_{s=1}^{Ns}$.
It is thus possible to see that (P1) with the empirical expectation~\eqref{eq:F-linear-appro}
can be equivalently re-written as
\vspace{.2cm}

\begin{subequations}\label{eq:AP1-all}
\hspace{-.5cm}
\fbox{
 \addtolength{\linewidth}{-2\fboxsep}%
 \addtolength{\linewidth}{-2\fboxrule}%
 \begin{minipage}{0.98\linewidth}
\begin{align}
\hspace{-0.3cm}\text{(AP1)}
\min_{\substack{\mathbf{p}_G,\mathbf{p}_W,\bm{\theta}, \\ \eta,\{\mathbf{v}_s,u_s\}_{s=1}^{Ns}}}
~& \sum_{m=1}^{M} C_m(p_{G_m}) + \mu\left(\eta + \frac{\sum_{s=1}^{Ns} u_s}{N_s(1-\beta)}\right) \label{eq:AP1-obj}\\
\st\quad\quad
& \eqref{eq:DCOPF-line}-\eqref{eq:DCOPF-pW} \nonumber \\
& \mathbf{v}_s \succeq \mathbf{p}_{W}-\mathbf{w}_s,~s = 1,\ldots,N_S \label{eq:AP1-v}\\
& \eta+u_s \ge  \mathbf{c}_W^\cT\mathbf{v}_s,~s = 1,\ldots,N_S \label{eq:AP1-u-eta} \\
& \mathbf{v}_s \succeq \mathbf{0}, u_s \ge 0,~s = 1,\ldots,N_S. \label{eq:AP1-uv}
\end{align}
\end{minipage}
}
\end{subequations}

\vspace{.2cm}

By introducing upper bounds $\{\mathbf{v}_s,u_s\}_{s=1}^{Ns}$, the non-smooth projection terms in the objective~\eqref{eq:DCOPF-obj} are
equivalently transformed to linear constraints~\eqref{eq:AP1-v}-\eqref{eq:AP1-uv}.
Thus, depending on whether $\{C_m(\cdot)\}_{m=1}^M$ are convex quadratic or piece-wise linear,
(AP1) is either a convex quadratic program (QP) or a linear program (LP),
which can be efficiently addressed by QP/LP solvers.

\begin{table}[t]
\centering
\caption{Generators Data.
The units of $\mathbf{p}_G^{\min(\max)}$, $c_m$ and $d_m$ are~MW, \$/(MWh)$^{2}$ and \$/MWh, respectively.}\label{tab:GEN}
\begin{tabular}{c||*{6}{c}}\hlinew{0.8pt}
Bus No.     &1      &2     &13         &22         &23         &27  \\   \hline
$\mathbf{p}_G^{\min}$     &0      &0     &0         &0         &0         &0 \\
$\mathbf{p}_G^{\max}$       &64    &64    &32        &40       &24    &44 \\
$c_m$                         &0.0200    &0.0175    &0.0250    &0.0625    &0.0250    &0.0083 \\
$d_m$                      &2.00     &1.75     &3.00     &1.00     &3.00     &3.25  \\
\hlinew{0.8pt}
\end{tabular}
\end{table}

\begin{table}[t]
\centering
\caption{Energy Purchase Prices and Forecast Wind Power.
The units of $c_{W_m}$ and $\bar{w}_m$ are~\$/MWh and MW, respectively. }\label{tab:W}
\begin{tabular}{l||*{7}{c}}\hlinew{0.8pt}
Bus No.          &1   &3   &7    &15     &19     &24     &26  \\ \hline
$c_{W_m}$      &3.5  &4.15    &2.65  &5.57   &4.64   &4.02   &6.75  \\
$\bar{w}_m$     & 6.00   & 0.31  & 7.66   & 8.01    &8.42   &8.44  & 8.46  \\
\hlinew{0.8pt}
\end{tabular}
\end{table}

Finally, it is worth pointing out that under mild conditions, the optimal solution set of (AP1) converges exponentially fast to its counterpart of (P1)
as the sample size $N_s$ increases. Due to space limitations, the proof of this claim is omitted. Interested readers are referred to~\cite{Kleywegt} for the
detailed analysis of the generic problem tackled using the theory of large deviations.

\section{Numerical Tests}\label{sec:Numericalresults}


Performance of the novel robust DC-OPF dispatch is corroborated via numerical tests using the IEEE $30$-bus benchmark system~\cite{PSTCA}.
The convex program (AP1) is solved using the~\texttt{CVX} package together with the~\texttt{SeDuMi} solver~\cite{cvx}, \cite{sedumi}.
The IEEE $30$-bus test system includes $41$ transmission lines and $6$ conventional generators the data of which are listed in Table~\ref{tab:GEN}.
The generation costs are $C_m(P_{G_m}) := c_m p_{G_m}^2 + d_m p_{G_m}$, for $m=1,\ldots,M$.
Other system parameters such as transmission line limits and base load demands are specified as in~\cite{MATPOWER}.

To simulate high penetration of wind energy, real data originally provided by Kaggle for the 
wind energy forecasting competition in $2012$ were utilized~\cite{kaggle}.
The dataset contains the hourly normalized power output of seven correlated wind farms.
They are assumed to be attached to different buses of the test system (cf.~Table~\ref{tab:W}).

Clearly, actual wind power output samples $\{\mathbf{w}_s\}_{s=1}^{N_s}$ in~\eqref{eq:AP1-v} are needed as the input of (AP1).
The required samples can be obtained via forecast wind power data, or, the distributions of wind speed together with the wind-speed-to-wind-power mappings
(cf.~\cite{YuNGGG-ISGT13}). To this end, the model $\mathbf{w}_s = \bar{\mathbf{w}} + \mathbf{n}_s$
is postulated to accomplish the sampling task. The day-ahead forecast wind power $\bar{\mathbf{w}} := [\bar{w}_1,\ldots,\bar{w}_M]^\cT$ is
chosen as the Kaggle data observed at $8$~A.M. of May $22, 2012$ (cf.~$\bar{w}_m$ in Table~\ref{tab:W}).
The forecast error $\mathbf{n}_s$ is assumed to be a zero-mean correlated Gaussian random vector for simplicity.
The covariance matrix of $\mathbf{n}_s$ was empirically estimated using Kaggle data across $589$ hours between $05/01/2012$ and $06/26/2012$.
Finally, negative-valued elements of the generated samples $\{\mathbf{w}_s\}_{s=1}^{N_s}$ were truncated to zero under physical constraints.
The probability level $\beta= 0.95$ and the sample size $N_s = 1,000$ were set in all the tests.

\begin{figure}[t]
\centering
\includegraphics[width=0.45\textwidth]{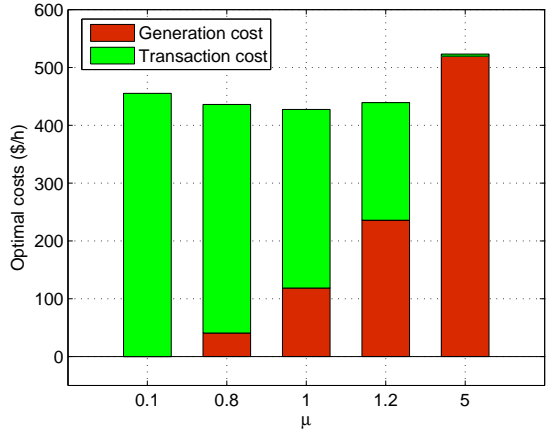}
\caption{Optimal costs for varying weight $\mu$.}
\label{fig:costAll}
\end{figure}
\begin{figure}[t]
\centering
\includegraphics[width=0.45\textwidth]{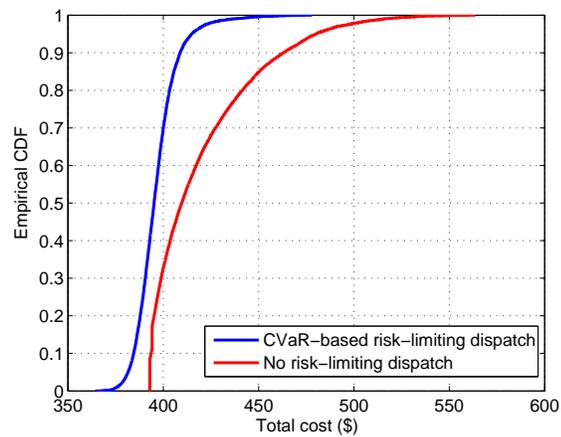}
\caption{Empirical CDFs of the optimal total costs.}
\label{fig:costCDF}
\end{figure}

\begin{table}[t]
\centering
\caption{Mean and variance of total costs: risk-limiting dispatch versus no risk-limiting dispatch.}\label{tab:meanVar}
\begin{tabular}{l||*{2}{c}}\hlinew{0.8pt}
Total cost         & Mean  & Variance \\
No risk-limiting dispatch         &419.87   & 856.24      \\
CVaR-based risk-limiting dispatch  & 396.40   &126.09   \\
\hlinew{0.8pt}
\end{tabular}
\end{table}

\begin{figure}[t]
\centering
\includegraphics[width=0.45\textwidth]{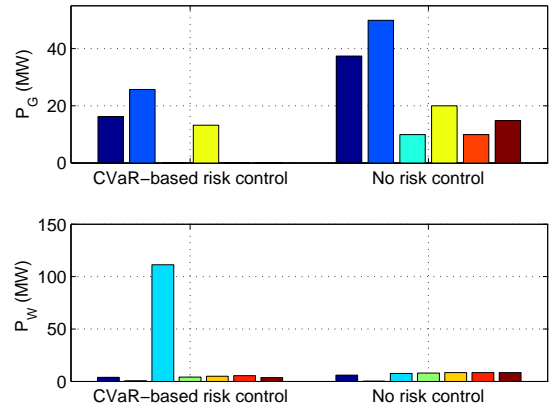}
\caption{Optimal power dispatch of $\mathbf{p}_{G}$ and $\mathbf{p}_{W}$.}
\label{fig:PGPW}
\end{figure}

\begin{figure}[t]
\centering
\includegraphics[width=0.45\textwidth]{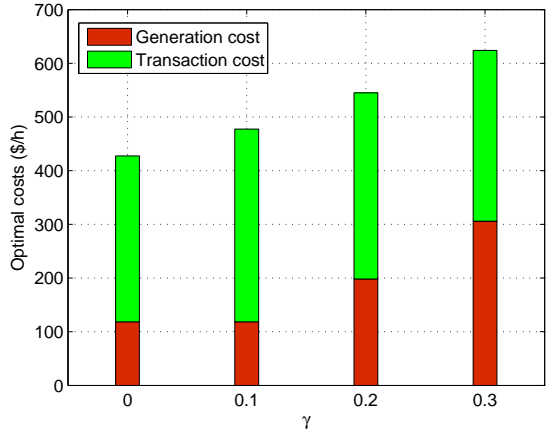}
\caption{Optimal costs for varying overload ratio $\gamma$.}
\label{fig:loadLift}
\end{figure}

\begin{figure}[t]
\centering
\includegraphics[width=0.45\textwidth]{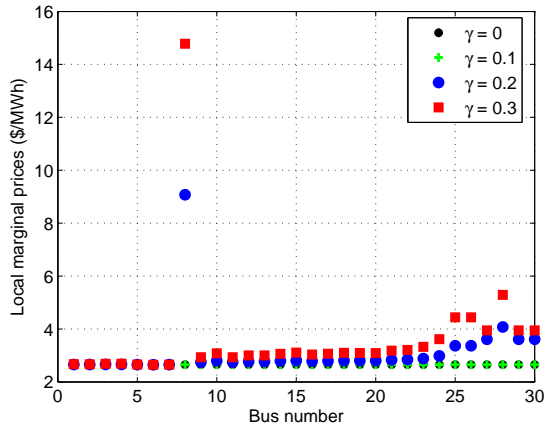}
\caption{Locational marginal prices for varying overload ratio $\gamma$.}
\label{fig:LMP-CVaR}
\end{figure}

Fig.~\ref{fig:costAll} depicts the optimal costs varying with the regularization weight $\mu$. 
It is clear that the conditional transaction cost decreases as $\mu$ increases. 
Since larger $\mu$ effects heavier penalty on the transaction cost, less $\mathbf{p}_{W}$ will be scheduled 
to reduce the risk of wind power shortage. This will result in the increase of the conventional generation 
$\mathbf{p}_{G}$ with the generation cost as shown in Fig.~\ref{fig:costAll}.

Fig.~\ref{fig:costCDF} shows the cumulative distribution functions (CDFs) of the optimal costs of 
the novel CVaR-based risk-limiting dispatch and no risk-limiting counterpart.
For the latter, the forecast wind power quantity $\bar{\mathbf{w}}$ is simply used in the 
nodal balance~\eqref{eq:DCOPF-node} for solving (P1), but without
the regularizer $F_{\beta}$. Note that after solving (P1) to obtain the optimal dispatch, 
the generation cost part becomes fixed. Hence, randomness of the total cost
comes from that of the transaction cost. Clearly, the no-risk control scheme always incurs 
a higher total cost than the novel CVaR-based risk-limiting approach. 
The values of the mean and variance of the optimal total costs are given in Table~\ref{tab:meanVar}, 
which again speak for the merits of the proposed scheme that exhibits reduced expected cost and variance.

Fig.~\ref{fig:PGPW} depicts the optimal power dispatch of conventional generation $\mathbf{p}_{G}$ and committed wind power $\mathbf{p}_{W}$ for both schemes:
CVaR-based risk control versus no-risk control. It can be seen that for the CVaR-based approach, the large scheduled wind power ${p}_{W_m}$ at bus $7$ makes the
$\mathbf{p}_{G}$ lower than that of the no-risk control, and thus gives rise to a lower optimal total cost as corroborated by the CDFs. This happens because
the proposed scheme can leverage the condition that purchase price $c_{W_m}$ at bus $7$ is the lowest one
among all wind power injection buses (cf.~Table~\ref{tab:W}), which allows for relatively high energy transaction at this bus.

Finally, the effects of demand overload are tested using the results of Figs.~\ref{fig:loadLift} and~\ref{fig:LMP-CVaR}.
Load demands at all buses are scaled up by $\gamma$, based on the original data of the IEEE $30$-bus system.
As expected, the total cost increases with the increase of the overload ratio $\gamma$ as confirmed by Fig.~\ref{fig:loadLift}.
It is interesting to observe that the overload hardly affects the transaction cost due to the CVaR-aware risk control mechanism.
Being important components of electricity markets, LMPs represent the cost (revenue) of buying (selling) electricity
at a particular bus~\cite{ExpConCanBook}. In Fig.~\ref{fig:LMP-CVaR}, LMPs are plotted for varying overload ratios $\gamma = 0, 0.1, 0.2, 0.3$.
Note that all LMPs should be equal in the case of no transmission line congestion.
Clearly, the congestion happens as overload demand increases.


\section{Conclusions and Future Work}\label{sec:Conclusion}
CVaR-based DC-OPF with wind integration was investigated in this paper. A convex optimization problem was formulated
considering the trade off between conventional generation cost and conditional energy transaction cost.
The CVaR-based regularizer plays an important role of risk-limiting dispatch, thus effecting smart utilization
of renewables to reduce the cost of conventional power generation, while taking limited risk of wind power shortage.

A number of appealing directions open up towards extending the proposed model and approach.
These include CVaR-based UC and AC-OPF, uncertain load demand and transaction costs, as well as distributed scheduling.

\section*{Acknowledgement}
The authors are grateful to Drs. Vassilis Kekatos and Nikolaos Gatsis for helpful discussions
and their aid with the data collection in~\cite{YuNKGG-DSP13}.
The authors also would like to thank Prof. Shuzhong Zhang of the Univ. of Minnesota,
for his inspiring suggestion to consider CVaR.

\bibliographystyle{IEEEtran}
\bibliography{biblio}
\end{document}